\documentclass[11pt]{article}
\usepackage{amssymb,amsfonts,amsthm}
\usepackage{color}
\usepackage[latin1]{inputenc}
\usepackage{makeidx}
\usepackage{amsmath}

\def\openC{{\rm C\kern-.18cm\vrule width.8pt height 7pt depth-.2pt \kern.18cm}}
\def\openN{{{\rm I}\kern-.16em {\rm N}}}
\def\openR{{{\rm I}\kern-.16em {\rm R}}}
\def\openT{{{\rm T}\kern-.42em {\rm T}}}
\def\openZ{{{\rm Z}\kern-.28em{\rm Z}}}

\def\eop{\hfill\rule{2.5mm}{2.5mm}}
\def\pf{\par\smallbreak\noindent {\bf Proof.} \ }
\newtheorem{thm}{Theorem}[section]
\newtheorem{cor}[thm]{Corollary}
\newtheorem{lem}[thm]{Lemma}
\newtheorem{prop}[thm]{Proposition}

\theoremstyle{definition}

\def\eop{\hfill\rule{2.5mm}{2.5mm}}

\begin{document}
\title{
{\textbf{\Large{On approximation tools and its applications on compact homogeneous spaces}}} \vspace{-4pt}
\author{A. O. Carrijo\thanks{E-mail: angelina.carrijo@usp.br}
\,\,\&\,\,
T. Jord\~{a}o\,\thanks{E-mail: tjordao@icmc.usp.br. 
Partially supported by FAPESP, grant $\#$ 2016/02847-9}}
}
\date{}
\maketitle \vspace{-30pt}
\bigskip

\begin{center}
\parbox{15 cm}{{\small \textsc{Abstract.} We prove a characterization for the Peetre type $K$-functional on $\mathbb{M}$, a compact two-point homogeneous space, in terms the rate of approximation of a family of multipliers operator defined to this purpose.\ This extends the well known results on the spherical setting.\
The characterization is employed to show that an abstract H\"{o}lder condition or finite order of differentiability condition imposed on kernels generating certain operators implies a sharp decay rates for their eigenvalues sequences.\ The latest is employed to obtain estimates for the Kolmogorov $n$-width of unit balls in Reproducing Kernel Hilbert Space (RKHS).}}
\end{center}


\thispagestyle{empty}

%
%
\section{Introduction}

The basic framework of this paper refers to a compact two-point homogeneous space of dimension $m\geq 1$.\ Denoting by $\mathbb{M}$ this space which is both a Riemannian $m$-manifold and a compact symmetric space of rank $1$ for which there is a well-developed harmonic analysis structure on them.\ A very large class of problems in approximation theory, harmonic analysis and functional analysis (as it can be seen in the present paper) can be considered naturally on these spaces.\ 

Each one of these manifolds $\mathbb{M}$ has an invariant Riemannian (geodesic) metric $d(\cdot,\cdot)$ which can be considered normalized so that all geodesics on $\mathbb{M}$ have the same length, namely, $2\pi$.\ Also $\mathbb{M}$ is endowed naturally with a measure $dx$ induced by the normalized left Haar measure which exists on a component of $\mathbb{M}$ seeing as a quotient.\ According to Wang \cite{wang}, compact two-point homogeneous spaces are: the unit spheres $\mathbb{S}^m$, $m=1,2,\dots$; the real projective spaces $\mathbb{P}^m(\mathbb{R})$, $m=2,3,\dots$; the complex projective spaces $\mathbb{P}^m(\mathbb{C})$, $m=4,6,\dots$; the quaternion projective spaces $\mathbb{P}^m(\mathbb{H})$, $m=8,12,\dots$ and $16$-dimensional Cayley's elliptic plane $\mathbb{P}^{16}$.\ These spaces have a very similar geometry between them and we shall assume here that $\mathbb{M}\neq S^m$ since the results we will present here already have their spherical version explored (see \cite{daidi, jordaosun} and references quoted there).

If we write $\mathcal{B}$ for {\em Laplace-Beltrami operator} on $\mathbb{M}$ it is known that its differential form depends on a pair of index $(\alpha,\beta)$ varying according to the space.\ Namely,
 $\alpha=(\sigma+\rho-1)/2=(m-2)/2$ and $\beta=(\rho-1)/2$ where: for $S^m$, $\sigma=0$, $\rho=m-1$; for $\mathbb{P}^m(\mathbb{R})$, $\sigma=m-1$, $\rho=0$; for  $\mathbb{P}^m(\mathbb{C})$, $\sigma=m-2$, $\rho=1$; for $\mathbb{P}^m(\mathbb{H})$, $\sigma=m-4$, $\rho=3$; and for $\mathbb{P}^{16}$, $\sigma=8$, $\rho=7$.\ We suggest \cite{kushpel2,platonov1,platonov2} and references therein for information above and details  about these spaces.\ 
 
Restricting ourselves to $m\geq2$ and $1\leq p\leq \infty$ we denote by $p'$ its exponent conjugate, i.e., $1/p+1/p'=1$.\ We write $(L^p(\mathbb{M}), \|\cdot\|_p)$ the usual Banach spaces of $p$-integrable complex functions on $\mathbb{M}$.\ In particular, $L^2(\mathbb{M})$ is a Hilbert space with the inner product $\langle\cdot ,\cdot \rangle_2$ defined by the normalized (by $\sigma_m$ the volume of $\mathbb{M}$) integral between square-integrable functions.

The Laplace-Beltrami operator on $\mathbb{M}$ has a discrete spectrum which arranged in an increasing order is given by $\{k(k+\alpha+\beta+1): k=0,1,\ldots\}$.\ For each $k$ the eigenspace $\mathcal{H}^m_k$ attached to $k(k+\alpha+\beta+1)$ has finite dimension $d_k^m:=\dim \mathcal{H}^m_k$ and they are mutually orthogonal.\ If we write $\{Y_{k,j}: j=1,2, \ldots, d_k^m\}$ for an orthonormal basis of $\mathcal{H}^m_k$, then $\{Y_{k,j}:  k=0,1,\ldots, \, j=1,2, \ldots d_k^m\}$ is an orthonormal basis of $L^2(\mathbb{M})$.\ This permits us to consider naturally Fourier expansions on $L^2(\mathbb{M})$.\ On the sphere all those objects are the well known space of spherical harmonics in $m+1$ variables and degree $k$ (\cite{dai-xu}).

The Fourier coefficients of a function $f \in L^p(\mathbb{M})$ are defined by
$$
\widehat{f}_{k,j}:= \frac{1}{\sigma_m}\int_{\mathbb{M}}f(y)\,\overline{Y_{k,j}(y)}\, dy, \quad j=1,2, \ldots, d_k^m,\quad  k=0,1,\ldots.
$$
We write $S_t(\cdot)$ for {\em shifting operator} (see \cite{dai, platonov1}) on $L^2(\mathbb{M})$, which is defined by the average of a function in a ``ring" of $\mathbb{M}$, namely for each $x\in \mathbb{M}$ the set is $\sigma_t^x:=\{y\in \mathbb{M}: d(x,y)=t\}$, $0<t<\pi$, with the induced measure.\ Then the addition formula also available in this context implies the following Fourier expansion of the shifting operator on $L^2(\mathbb{M})$:
\begin{equation}\label{shiftingfourier}
S_t(f) \sim \sum_{k=0}^{\infty} Q^{(\alpha,\beta)}_k(\cos t)\, \mathcal{Y}_k(f), \quad  f\in L^2(\mathbb{M}),
\end{equation}
where $Q^{(\alpha,\beta)}_k$ denotes the normalized Jacobi polynomial, it means $Q^{(\alpha,\beta)}_k(1)=1$, and $\mathcal{Y}_k$ is the projection of $L^2(\mathbb{M})$ onto $\mathcal{H}^m_k$, $k=0,1,\ldots$.\ All the facts mentioned above can be found constructed and/or explored in the cited references and \cite{kushpel2, platonov2}. 

We write $\mathcal{B}^{r}(f)$ to denote the \emph{fractional derivative of order $r$} which is defined on $\mathbb{M}$ in the distributional sense and given by $\mathcal{B}^{r}(f)\sim\sum_{k=0}^{\infty}(k(k+\alpha+\beta+1))^{r/2}\,\mathcal{Y}_k(f)$, we are allowed to consider the Sobolev class
$$
W_p^r(\mathbb{M}):=\left\{f\in L^p(\mathbb{M}):\mathcal{B}^r(f)\in L^p(\mathbb{M})\right\},
$$ 
endowed which the usual norm $\|\cdot\|_{W_p^r}:=\|\cdot\|_p+\|\mathcal{B}^r(\cdot)\|_p$.

Let us consider $r>0$,  $t>0$ and $f\in L^p(\mathbb{M})$.\ We introduce the Peetre-type {\it $K$-functional of fractional order $r$} given by
\begin{equation}\label{Kfunctional}
K_r(f,t)_p:= \inf_{g\in W_p^r(\mathbb{M})}\left\{\|f-g\|_p+t^r\|g\|_{W_p^r}\right\}.
\end{equation}
The \textit{$r$-th moduli of smoothness} 
\begin{equation}\label{moduli}
\omega_{r}(f,t)_p:=\sup\left\{\|(I-S_s)^{r/2}(f)\|_p: s \in (0,t]\right\}.
\end{equation}
And the \textit{generalized shifting operator}
\begin{equation}\label{combshift}
S_{r,t}(f):= \frac{-2}{{2r \choose r}}\sum_{j=1}^{\infty}(-1)^j{2r \choose r-j}S_{jt}(f),
\end{equation}
where for $r,s$ real numbers
$$
{r \choose s} = \frac{\Gamma(r+1)}{\Gamma(s+1)\Gamma(r-s+1)}, \,\,\,\, \mbox{for $s\not\in \mathbb{Z}_{-}$,} \,\,\,\, {r \choose 0}=r \,\,\,\, \mbox{and} \,\,\,\, {r \choose s}=0, \,\,\,\, \mbox{for $s\in \mathbb{Z}_{-}$.}
$$
It is not difficult to see that this operator well defined and a bounded operator on $L^p(\mathbb{M})$, $1\leq p\leq \infty$.\ 

Platonov showed that the $K$-functional of fractional order and the moduli of smoothness are related in a asymptotic sense.\ Notation $A(t)\asymp B(t)$ stands for $B(t)\lesssim A(t)$ and $A(t)\lesssim B(t)$, while $A(t)\lesssim B(t)$ means that $A(t)\leq c\, B(t)$, for some constant $c\geq 0$ not depending upon $t$.

\begin{thm}(\cite[Theorem 1.2]{platonov1}) For $1<p<\infty$ and $r\geq 1$ a natural number, it holds\label{equivalences}
$$
K_{2r}(f,t)_p\asymp \omega_{2r}(f,t)_p \quad f\in L^p(\mathbb{M}),\,  t>0.
$$
\end{thm}

Our main interest on these tools is its relation with the decay of Fourier coefficients of functions in terms of the rate of approximation of generalized shifting operator.\ It has shown extremely an important and efficient tool to get good estimates for both Fourier coefficients of functions satisfying a generalized H\"{o}lder condition and eigenvalues sequences of positive integral integral operators with H\"{o}lderian kernels (see \cite{jordaosun}).\ The relation we have stablished is the following.

\begin{thm}\label{desejada} For $1<p<\infty$ and $r\geq 1$ a natural number, it holds\label{equivalences}
$$
K_{2r}(f,t)_p\asymp \|S_{r,t}(f)-f\|_p, \quad f\in L^p(\mathbb{M}),\,  t>0.
$$
\end{thm}
This extends the spherical version of it which can be found in \cite{daidi}.\ A simple adaptation of well-known results on spherical functions, namely the Hausdorff-Young inequality, on $\mathbb{M}$ play an important role to show the following decay of Fourier coefficients.\ He fix the notation $
s_k(f):=\sum_{j=1}^{d_k^m}|\widehat{f}_{k,j}|^2$, $k=0,1,\ldots$.

\begin{thm} \label{ournequality} Let $r$ be a positive interger.\ 
If $p \in (1,2]$, then
$$
\left\{\sum_{k=1}^{\infty}(d_k^m)^{(2-p')/2p'}(\min\{1, tk\})^{rq} [s_k(f)]^{p'/2}\right\}^{1/q}\lesssim \|S_{r,t}(f)-f\|_p,\quad f \in L^p(\mathbb{M}),
$$
The inequality above becomes an equality in the case $p=2$.\ And,  if $p=1$, then
$$
\sup_{k\geq 0}\left\{(d_k^m)^{-1/2}(\min\{1, tk\})^{rq}[s_k(f)]^{1/2}\right\}\lesssim \|S_{r,t}(f)-f\|_1,\quad f \in L^1(\mathbb{M}).
$$
\end{thm}

Results above permit us to analyze the asymptotic behavior of eigenvalues sequences
 $\{\lambda_n\}_n$ of certain integral operators generated by H\"{o}lderian kernels having a Mercer-like series expansion.\ For a historical review on the spherical setting we suggest see \cite{jordaosun}.

We will be dealing with integral operators 
$\mathcal{L}_K(f)=\int_{\mathbb{M}} K(\cdot,
y)f(y)\,dy$,  having the generating kernel $K\colon
\mathbb{M}\times \mathbb{M}\longrightarrow\mathbb{C}$ belonging to
$L^2(\mathbb{M}\times \mathbb{M})$.\ It is easily seen that $\mathcal{L}_K$ defines a compact operator on $L^2(\mathbb{M})$.\ The study concerns to kernels on $\mathbb{M}\times \mathbb{M}$ of the form:
\begin{equation}\label{kernel}
K(x, y) = \sum_{k=0}^\infty \sum^{d_k^m}_{j=1} a_{k,j}\,  Y_{k,j}(x)\, Y_{k,j}(y), \quad \sum_{k=0}^\infty \sum^{d_k^m}_{j=1} a_{k,j}<\infty, \,\,\, x,y\in \mathbb{M}.
\end{equation}
We work under two basic conditions: the first one, called \emph{positivity}, means that the expansion coefficients are non-negative, i.e., $a_{k,j} \geq 0$; and the second one, called \emph{monotonicity}, means that the expansion coefficients are monotone decreasing with respect to $k$, i.e., $a_{k+1,j} \le a_{k,j'},\;
1 \le j, j' \le d_k^m$.\  

Recently, Berg and collaborators (\cite{berg}) showed the characterization given by Schoenberg (\cite{schoenberg}) for continuous zonal positive definite kernels on the sphere as series expansion given by formula (\ref{kernel}) with coefficients do not depending on index $j$ and satisfying the positivity definition above for positive definite kernels, holds in a general setting, namely on products of compact Gelfand pairs with locally compact groups.\ Therefore, assumptions made here on compact two-point homogeneous spaces are very natural and expected in most of the applications.  

The first application is continuation of the designed in \cite{jordaosun}.\ We say that a kernel $K$ on $\mathbb{M}$ satisfies the $(B, \beta)-$H\"{o}lder condition if there exist
a $\beta \in (0,2]$ and a function $B$ in $L^1(\mathbb{M})$ such that
\begin{equation}\label{SScondition}
\left|S_t(K(y,\cdot))(x)-K(y,x)\right| \leq B(y)\,t^{\beta},\quad x,y \in
\mathbb{M}, \,\, t\in (0,\pi).
\end{equation}
It is not hard to see that this definition is a generalized version the usual H\"{o}lder condition.

Assumption of positivity on $K$ implies self-adjointness of the integral operator, then the standard spectral theorem is applicable and we obtain a sequence of nonnegative real numbers $\{\lambda_n(\mathcal{L}_K)\}_n$ which is the eigenvalues sequence of $\mathcal{L}_K$.
 
\begin{thm}\label{thmdecay2} Let ${\cal{L}}_K$ be the integral operator induced by a kernel $K$ as in (\ref{kernel}) and under assumptions of positivity and monotonicity.\ If $K$ satisfies the $(B,\beta)$-H\"{o}lder condition, then it holds
$$
\lambda_n({\cal{L}}_K)=O( n^{-1-\beta/m}), \quad \mbox{as} \,\,\, n\to \infty.
$$
\end{thm}

For a positive real number $r$, we write $\mathcal{B}^{r,0}K$ for  the action of the fractional derivative operator only applied to the first variable.\ Also, $K^y$ denotes the function $ \cdot \mapsto K(\cdot, y)$, $y\in \mathbb{M}$.\ We warn the reader that the terminology ``trace-class", common in operators theory appears in the next corollary, it basically means that the trace of the operator is finite and independent of the choice of basis.\ See \cite{mathcomp}, for details and applications.

\begin{cor}\label{thmdecay1} Let ${\cal{L}}_K$ be the integral operator induced by a kernel $K$ as in (\ref{kernel}) and under assumptions of positivity, monotonicity and such that for a fixed $r>0$, all $K^y$ belong to $W_2^{2r}(\mathbb{M})$.\ If the integral operator generated by $\mathcal{B}^{2r,0}K$ is trace-class, then
$$
\lambda_n({\cal{L}}_K)=O( n^{-1-2r/m}), \quad \mbox{as} \,\,\, n\to \infty.
$$
\end{cor}

This extends both Theorem 2.5 in \cite{mathcomp} and Theorem 3 in \cite{jordaosun} for compact two-point homogeneous spaces.\ We present it as consequence of previous theorem not because it is an immediate  consequence of it, but we apply similar techniques in order to prove it, though. 
 
The paper is organized as following.\ Section 2 contains a short description of the Hausdorff-Young type inequality its implication on the relation between the decay of Fourier coefficients of a function and the proof of Theorem \ref{ournequality}.\ Several technical lemmas are proved in order to present the proof of Theorem \ref{desejada}.\ A technique involving relations between the decay of Fourier coefficients and eigenvalues sequences of the operator is employed to prove \ref{thmdecay2} and \ref{thmdecay1}.\ In Section 4 we give a shortly background for Kolmogorov $n$-widths and apply our achievements to get sharp estimates for the Kolmogorov $n$-width of RKHS of H\"{o}lderian kernels.\ Finally, an example is given.

\section{Decay of Fourier coefficients}
\setcounter{equation}{0}

In this section we  present some background material in order to prove Theorema \ref{ournequality} and \ref{desejada}.\ They include realization theorem, moduli of smoothness and the associated $K$-functional as well.\ Relations between these were proved recently by Dai, Ditzian and Platonov on two-point homogenous spaces and play an important role here.\ References are \cite{dai,ditzian1,platonov2, platonov1}.\ 

A linear operator $T$ on $L^p(\mathbb{M})$ is called a {\em multiplier operator} if there exists a sequence $\{\mu_k\}_k$ of complex numbers such that ${\cal Y}_k(Tf)=\mu_k \, {\cal Y}_k (f)$, $k=0,1,\ldots$, for any $f \in L^p(\mathbb{M})$ and $T$ is bounded.\ In this case the sequence $\{\mu_k\}_k$ is called the sequence of multipliers of $T$.

An important property involving the $K$-functional is the Realization Theorem for $K_r(f,t)_p$ (\cite{ditzian1}), which is given by the relation below.\ In its statement, the multiplier operator $\eta_t$ depends upon a best approximation function $\eta\in C^{\infty}[0,\infty)$ such that $\eta=1$ in $[0,1]$, $\eta=0$ in $[2,\infty)$ and $\eta(s)\leq 1$, $s\in (1,2)$.\ The operator $\eta_t$ is defined by the formula
$\eta_t(f)=\sum_{k=1}^{\infty}\eta(tk)\,\mathcal{Y}_k(f)$ for all $f \in L^p(\mathbb{M})$.\ For $r>0$ and $f \in L^p(\mathbb{M})$ Realization Theorem (\cite{ditzian1}) assures that the $K$-functional $K_r(f,t)_p$ assumes its infimum via the operator $\eta_t$ as bellow:
\begin{equation}\label{realizationthm}
\|f-\eta_t(f)\|_p+t^r\,\|\eta_t(f)\|_{W_p^r}\asymp K_r(f,t)_p, \quad t>0.
\end{equation}

\begin{lem} \label{diti}(Hausdorff-Young type inequality) Let $q$ be the conjugate exponent of $p$.\  Then
$$
\left\{\sum_{k=1}^{\infty}(d_k^m)^{(2-q)/2q}\,[s_k(f)]^{q/2}\right\}^{1/q}\lesssim \|f\|_p,\quad 1<p\leq 2, \,\, f\in L^p(\mathbb{M}); 
$$
and
$$
\sup_{k\geq 0}\left\{(d_k^m)^{-1/2}[s_k(f)]^{1/2}\right\}\lesssim \|f\|_1,\quad f\in L^1(\mathbb{M}).
$$
\end{lem}

The proof is based on the Riez-Thorin interpolation and nothing different from the spherical setting (see \cite{ditzian2}, for example) that is why it is omitted here.

The following theorem relates the decay of the Fourier coefficients of a function to the rate of approximation of operator defined in formula (\ref{combshift}).\ In \cite{jordaomen} a proof of similar result is presented for a multiplier operator on the spherical setting and it is slightly different from below.

\vspace{0.25cm}

\noindent{\bf Proof of Theorem \ref{ournequality}.} Formula (\ref{shiftingfourier}) implies
\begin{eqnarray*}
S_{r,t} = \sum_{k=0}^{\infty}\left[\frac{-2}{{2r \choose r}}\sum_{j=1}^{r}(-1)^j{2r \choose r-j}Q^{(\alpha,\beta)}_k(\cos (jt))\right] \,\mathcal{Y}_k, \quad t >0.
\end{eqnarray*}
It means that $S_{r,t}$ is a multiplier operator and its multiplier sequence $\{m_r(k,t)\}_k$ is given by 
\begin{equation}\label{multiplier}
m_r(k,t) := \frac{-2}{{2r \choose r}}\sum_{j=1}^{r}(-1)^j{2r \choose r-j}Q^{(\alpha,\beta)}_k(\cos (jt)), \quad k=0,1\ldots, \,\, t>0.
\end{equation}

Fixing $f\in L^p(\mathbb{M})$, the linearity of the orthogonal projections imply that
$$
{\cal Y}_k(S_{r,t}(f)-f) =(m_r(k,t)-1){\cal Y}_k(f),  \quad k=0,1\ldots,
$$
whence
$$
\sum_{j=1}^{d_k^m}\widehat{(S_{r,t}(f)-f)}_{k,j}Y_{k,j}=(m_r(k,t)-1)\sum_{j=1}^{d_k^m}\hat{f}_{k,j} Y_{k,j},\quad k=0,1\ldots.
$$
Computing the $L^2$-norm of both sides we obtain $s_k(S_{r,t}(f)-f)=(m_r(k,t)-1)^2 s_k(f)$, $k=0,1\ldots,$
that is,
$$
(d_k^m)^{(2-q)/2q}\left[s_k(S_{r,t}(f)-f)\right]^{q/2}=(d_k^m)^{(2-q)/2q}|m_r(k,t)-1|^q [s_k(f)]^{q/2}.
$$
Taking in account that $1-m_r(k,t)\asymp (\min\{1, tk\})^r$ (this important and nontrivial equivalence is obtained in Lemma \ref{eqmultiplier}) we have
\begin{eqnarray*}
(d_k^m)^{(2-q)/2q}\left[s_k(S_{r,t}(f)-f)\right]^{q/2} &= & (d_k^m)^{(2-q)/2q}|m_r(k,t)-1|^q [s_k(f)]^{q/2}\\ & \asymp & (d_k^m)^{(2-q)/2q}(\min\{1, tk\})^{rq} [s_k(f)]^{q/2}.
\end{eqnarray*}
Finally Hausdorff-Young type formula reach us to the first inequality in the statement of the theorem.\ As for the equality assumption in the case $p=2$, it suffices to apply Parseval's identity in the equality above.\ The inequality in the case $p=1$ is settled in a similar fashion.\ \eop

Theorem above and Theorem \ref{desejada} permit us to choose the more convenient tool in order to study the decay of Fourier coefficients.\ We can also relate the decay of the Fourier coefficients of a function to the $K$-functional defined in (\ref{Kfunctional}).\ Ditzian \cite{ditzian1} proved this theorem on the spherical setting and for the special case in which $r$ is a positive integer.\ He remarks that the same proof can be slightly modified to fit for $r$ been a real number.\ Since the proof is a direct application of Hausdorff-Young type inequality (Lemma \ref{diti}) we choose do not reproduce it here (it can be founded in \cite[p. 9]{jordaosun}on the spherical setting).\ 

\begin{prop}\label{ditinequality1}  If $f$ belongs to $L^p(\mathbb{M})$, $1\leq p\leq 2$, $q$ is the conjugate exponent of $p$ and $r > 0$, then
\begin{equation}
\left\{\sum_{k=1}^{\infty}(d_k^m)^{(2-q)/2q}\,(\min\{1, tk\})^{rq}\,[s_k(f)]^{q/2}\right\}^{1/q}\lesssim K_r(f,t)_p,\quad t>0.
\end{equation}
\end{prop}

\subsection{Proof of Theorem \ref{desejada}}

The technic employed here is to get good estimates for the multiplier sequence attached to averaged operator and through an application of the Marcinkiewicz's Multiplier Theorem we prove Theorem \ref{desejada}.\ This proof is highly technical and we present it by steps.\ 

The first technical lemma brings estimates for the difference operator applied to Jacobi polynomials.\ The difference operator is defined inductively as follows: for a sequence $\{b_k\}_k$, we set $\triangle^0b_k=b_k$ and $\triangle b_k=b_{k+1}-b_k$ from that if $j$ is a positive integer $\triangle^jb_k=\triangle(\triangle^{j-1}b_k)$.

\begin{lem}\label{lem1}
If $\alpha\geq\beta\geq -\frac{1}{2}$, then 
\begin{eqnarray*}
\left|\triangle^jQ^{(\alpha,\beta)}_k(\cos t)\right| \lesssim \left \{ \begin{matrix} t^j,  & & kt \leq  1 \\ 
(kt)^{-(\alpha+1/2)},   & & kt \geq  1. \end{matrix} \right.
\end{eqnarray*}
\end{lem}

Its proof follows directly from \cite[Lemma 2]{dai-wang}.\ We apply it in order to obtain related estimates for the sequence $\{m_r(k,t)\}_k$ as follows.

\begin{lem}\label{lem2}
Let $\{m_r(k,t)\}_k$ be the multiplier sequence of operator $S_{r,t}$.\ If $0<t\leq \frac{\pi}{2r}$ and $j$ is a positive integer, then it holds
\begin{eqnarray*}
\left|\triangle^j m_r(k,t)\right| \lesssim \left \{ \begin{matrix} t^j, & 0<kt\leq 1 \\ t^j(kt)^{-(\alpha+1/2)}, & kt\geq 1. \end{matrix} \right.
\end{eqnarray*}
\end{lem}

\pf For each $k$ the representation of $m_r(k,t)$, given in formula (\ref{multiplier}), implies the following inequality
\begin{equation*}
\left|\triangle^j m_r(k,t)\right| \leq \frac{2}{{2r \choose r}}\sum_{j=1}^{r}\left|{2r \choose r-j}\right|\left|\triangle^j  Q^{(\alpha,\beta)}_k(\cos (jt))\right|.
\end{equation*}
An application of Lemma \ref{lem1} is enough to complete the proof. \eop

The next result asserts that we can represent the normalized Jacobi polynomial by a sum of cosines with nonnegative coefficients (we warn the reader that it does not hold for all Jacobi polynomial, see \cite{askey} for details).\ The proof for this fact can be found in \cite[p. 63--66]{askey} and we just need to observe that the one-dimensional unit sphere $S^1$ is isometrically embedded in any $\mathbb{M}$.\ 

\begin{lem}\label{represent} If $\alpha$ and $\beta$ are as described Section 1, then
$$
Q^{(\alpha,\beta)}_k(\cos \theta) = \sum_{v=0}^{k}\sum_{i=0}^{[v/2]} a_{v,i}\,(\cos (v-2i)\theta), \quad k=0,1,\ldots,
$$
where $a_{v,i}\geq 0$, $v=0,1,\ldots, k$ and $i=0,1,\ldots, v$.
\end{lem}

Explicitly, a simple calculation for $v=0,1,\ldots, k$ and $i=0,1,\ldots, v$, implies 
$$
a_{v,i}= b_{v,i}\frac{\Gamma(k+\alpha+1)\Gamma(\alpha-\beta+1)\Gamma(2\alpha+v+1)\Gamma(2\alpha+2v+2)\Gamma(k+\alpha+\beta+v+1)}{\Gamma(\alpha-\beta-k+v+1)\Gamma(k-v+1)\Gamma(\alpha+v+1)\Gamma(2\alpha+2v+1)\Gamma(k+2\alpha+v+2)},
$$
where $\Gamma (\cdot)$ stands the Gamma function and $b_{v,i}=c_{v,i}P_v^{(\alpha,\alpha)}(1)/P_k^{(\alpha,\beta)}(1)$ and $c_{v,i}$ are given by the Gegenbauer polynomial with index $\alpha$ representation (\cite[p. 93]{szego}) in terms of cosine.

The main idea behind the next result is that if $\{m_r(k,t)\}_k$ is the sequence of multipliers of $S_{r,t}$, then $1-m_r(k,t)\asymp (\min\{1, tk\})^r$.

\begin{lem}\label{eqmultiplier}
For $t\in[0,\pi/2]$ it holds
\begin{equation}\label{1}
0<a\leq \frac{1-m_r(k,t)}{(kt)^{2r}}\leq b<\infty, \quad \mbox{ for }\,\,  0<kt\leq \pi,
\end{equation}
where $a$ and $b$ are constants.\ Additionally, for any $\tau>0$ there exists $v_{r,\tau}<1$ such that
\begin{equation}\label{2}
m_r(k,t) \leq v_{r,\tau}, \quad \mbox{ for }\,\,  kt\geq \tau >0.
\end{equation}
\end{lem}

\pf  The main ideia of the proof is borrowed from \cite[Lemma 4.4]{daidi} but several considerations are needed.\ By Lemma \ref{represent} we have
\begin{eqnarray*}
1-m_r(k,t) &=& 1+ \frac{2}{{2r \choose r}}\sum_{j=1}^{r}(-1)^j{2r \choose r-j}\sum_{v=0}^{k}\sum_{i=0}^{[v/2]} a_{v,i}\,(\cos (v-2i)jt).
\end{eqnarray*}
From above the representation of a power of the sine function in terms of cosine implies 
\begin{eqnarray*}
1-m_r(k,t) 	&=&\frac{4^r}{{2r \choose r}}\sum_{v=0}^{k}\sum_{i=0}^{[v/2]} a_{v,i}\left[\sin^{2r}\left((v-2i)\frac{t}{2}\right)\right].
\end{eqnarray*}

If $kt\leq \pi$, then $\sin^{2r}\left(\frac{(v-2i)t}{2}\right)\leq \left(kt/2\right)^{2r}$.\ Which implies 
\begin{eqnarray*}
1-m_r(k,t) 	
\leq \frac{4^r}{{2r \choose r}}\sum_{v=0}^{k}\sum_{i=0}^{[v/2]} a_{v,i} \left(\frac{kt}{2}\right)^{2r}  =   \frac{(kt)^{2r}}{{2r \choose r}},
\end{eqnarray*}
and the proof of the right-hand side of inequality (\ref{1}) follows.\ 

On the other hand, starting from
\begin{eqnarray*}
1-m_r(k,t)= \frac{4^r}{{2r \choose r}}\sum_{v=0}^{k}\sum_{i=0}^{[v/2]} a_{v,i}\left[\sin^{2r}\left((v-2i)\frac{t}{2}\right)\right],
\end{eqnarray*}
for each index $i\leq [v/4]$ it holds $\sin^{2r}\left(\frac{(v-2i)t}{2}\right) \geq \left(vt/2\pi\right)^{2r}$. 
Then, we have for some positive constant $c$ that
\begin{eqnarray*}
1-m_r(k,t) \geq  \frac{4^r}{{2r \choose r}}\sum_{v=1}^{k}\sum_{i=0}^{[v/4]} a_{v,i}\,\left(\frac{vt}{2\pi}\right)^{2r}\geq  c\,  \frac{4^r}{{2r \choose r}}\left(\frac{kt}{2\pi}\right)^{2r}.
\end{eqnarray*}
Which finishes the proof of inequalities in formula (\ref{1}) in the statement of the lemma.

Let any $0<\tau\leq \pi$ and $kt \leq \pi$, inequality (\ref{1}) implies
$$
0<a \leq \frac{1-m_r(k,t)}{(kt)^{2r}} \leq \frac{1-m_r(k,t)}{\tau^{2r}}.
$$
And, therefore $m_r(k,t) \leq 1 - a\tau^{2r}$, if choose $v_{r,\tau}= 1-a\tau^{2r} <1$, it is (\ref{2})for $0<\tau \leq kt$.

Now if $kt\geq \pi$ we have
\begin{eqnarray*}
1-m_r(k,t) 	&\geq & \frac{4^r}{{2r \choose r}}\sum_{v=0}^{k}\sum_{i \in I(k)} a_{v,i}\left[\sin^{2r}\left((v-2i)\frac{t}{2}\right)\right],
\end{eqnarray*}
where 
$$
I(k) := \bigcup_{l=0}^{\left[\frac{vt}{2\pi}-\frac{1}{2}\right]}\left\{i: 0\leq i\leq [v/2];\quad \frac{\pi}{4}+l\pi \leq (v-2i)\frac{t}{2} \leq \frac{3\pi}{4}+l\pi \right\}.
$$
Writing $c'=\sum_{v=0}^{k}\sum_{i \in I(k)} a_{v,i}$ we have $1-m_r(k,t) \geq  \frac{4^r}{{2r \choose r}}\,c'$
And the lemma is proved.\ \eop

Now we present some more properties of the difference operator. 
\begin{lem}\label{lem7}\label{inducao} Let $\{a_k\}_k$, $\{b_k\}_k$ sequences of real numbers and $j$ a positive integer.

\noindent {\bf a)} It holds
\begin{eqnarray*}
\triangle^j(a_k b_k) = \sum_{i=0}^j{j \choose i}(\triangle^{j-i}a_k)(\triangle^ib_{k+j-i}).
\end{eqnarray*}

\noindent{\bf b)} If the sequence $\{a_k\}_k$ satisfies $a_k\geq a>0$, $k=0,1,\ldots$, then
\begin{eqnarray*}
|\triangle^{j}a_k^{-1}| \leq \frac{1}{|a_k|} \sum_{i=0}^{j-1} {j \choose i} |\triangle^i a_k^{-1}|\,|\triangle^{j+i}a_{k+i}|
   \leq  c \, \max_{0\leq i\leq j}  |\triangle^i a_k^{-1}|\,|\triangle^{j+i}a_{k+i}|,
\end{eqnarray*}
where $c= 2^j/a$.
\end{lem}

The proof of item a) follows by mathematical induction.\ For the part b) note that $a_k^{-1}a_k=1$, choosing the sequences $\{a_k^{-1}\}_k$ and $\{a_k\}_k$ we apply item a) and the proof follows.

\begin{lem}\label{lem9}
If $t\in[0,\pi/2]$, then for any positive integer $j$ and $\tau >0$ such that $0<kt<\tau$ the following holds
\begin{equation}\label{eq-lem8}
\left| \triangle^j \frac{1-m_r(k,t)}{(k(k+\alpha+\beta+1)t^2)^r} \right|  \lesssim 
\left[k^{-j}+k^{-j-1}\right].
\end{equation}
\end{lem}

The proof is omitted  since it can found in \cite[p. 255-258]{dai-xu}.\ We observe that inequality (\ref{eq-lem8}) holds for the sequence of multiplicative inverse sated as well.

The very last result we need in oder to present the proof of Theorem \ref{desejada} is the Marcinkiewicz's Multiplier Theorem for the compact two-point homogeneous spaces.\ This result gives us a sufficient condition such that a given operator constructed via sequences (multipliers) be bounded.

\begin{thm}{\em \cite[Theorem 7.1]{Bonami}} \label{Marcinkiewicz}
Let $\mathbb{M}$ be a compact two-point homogeneous space of dimension $m$ and $\{\mu_j\}_j$ a 
sequence of real numbers satisfying

\medskip
\noindent {\bf i)} $\sup_j\{|\mu_j|\}\leq M<\infty$;

\medskip
\noindent {\bf ii)} $\sup_j\left\{2^{j(s-1)}\sum_{l=2^j}^{2^{j+1}}|\bigtriangleup^s\mu_l|\right\}\leq M<\infty$,
with $s=(m+1)/2$ if $s$ is odd and $s=(m+2)/2$ if $s$ is even.\ Then, it holds 
\begin{eqnarray*} 
\left\| \sum_{k=0}^{\infty} \mu_k \, \mathcal{Y}_k(f) \right\|_p \leq a_p\,M\,\|f\|_p, \quad f \in L^p(\mathbb{M}),\,\, 1<p<\infty,
\end{eqnarray*}
where $a_p$ is a constant which does not depend on $f$.\ 
\end{thm}
\vspace{0.25cm}

\noindent{\bf Proof of Theorem \ref{desejada}.} 
In order to prove the equivalence stated, it is enough to show that for some positive $a$ the following three inequalities hold
\begin{equation}\label{eq1}
 \| f- \eta_{at} f \|_p \lesssim  \| f- S_{r,t} f \|_p,
\end{equation}
\begin{equation}\label{eq2}
t^{2r}\, \| \mathcal{B}^{2r} (\eta_{at} f) \|_p \lesssim \| f- S_{r,t} f \|_p
\end{equation}
and
\begin{equation}\label{eq3}
\| \eta_{at} f - S_{r,t} (\eta_{at} f) \|_p \lesssim t^{2r}\,  \| \mathcal{B}^{2r} (\eta_{at} f) \|_p.
\end{equation}
Also, observe that (\ref{eq1}) is assured if
\begin{equation} \label{eq4}
\| f - \eta_{at} f - ( I + S_{r,t} + \cdots + S_{r,t}^4 ) \,(I- \eta_{at})\, (f- S_{r,t} f )\|_p \lesssim \| f- S_{r,t} f \|_p.
\end{equation}

To obtain inequality above we need to show that 
\begin{equation*}
\mu_{k,1}: = \left(1-\eta(atk)\right)\,\frac{m_r(k,t)^4}{1-m_r(k,t)},\quad k=0,1,\ldots,
\end{equation*}
define sequence satisfying conditions in Theorem \ref{Marcinkiewicz}.\ We first note that if $atk\leq 1$, $\eta(atk)=1$ and then, $\mu_{k,1}=0$.\ If $atk\geq \tau >1$, from Lemma \ref{eqmultiplier} we have $1-m_r(k,t)\geq c_{\tau,r}>0$, also it is clear that $|1-\eta(atk)|\leq c$ for some constant $c$ and since $\{(1/{kt})^{4\alpha+4/2-s}\}_k$ is bounded for $4\alpha+4/2-s\geq 0$, Lemma \ref{lem2} implies 
\begin{eqnarray*}
|\triangle^s \mu_{k,1}|  \lesssim   |\triangle^s m_r(k,t)^4|  \lesssim \frac{t^s}{(kt)^{4\alpha+4/2}} \lesssim \left(\frac{1}{kt}\right)^{4\alpha+4/2-s} \frac{1}{k^s}  \lesssim \frac{1}{k^s}.
\end{eqnarray*}
From the inequality above, we have
\begin{eqnarray*}
\sup_j \left\{2^{j(s-1)} \sum_{k=2^j}^{2^{j+1}} |\bigtriangleup^s \mu_{k,1}|\right\} \lesssim  \sup_j \left\{2^{j(s-1)} \sum_{k=2^j}^{2^{j+1}}   \frac{1}{k^s}\right\} 
 			& \leq&  \sup_j 2^{j(s-1)} \sum_{k=2^j}^{2^{j+1}} \frac{1}{(2^j)^s}\lesssim 1.
\end{eqnarray*}
And therefore Theorem \ref{Marcinkiewicz} assures that $\{\mu_{k,1}\}_k $ is a multiplier sequence and inequality (\ref{eq1}) is proved.

Heading to inequality (\ref{eq3}), we proceed analogously as above but taking in account the convenient multiplier sequence.\ For $0<kt<\tau$, $\tau>0$ let 
\begin{equation*}
\mu_{k,2} := \frac{1-m_r(k,t)}{(k(k+\alpha+\beta+1))^rt^{2r}} \, \eta(atk), \quad k=1,2,\ldots
\end{equation*}
be the sequence of  multipliers to application involved in inequality (\ref{eq3}).\ 
We write 
$$
a_k:=  \frac{1-m_r(k,t)}{(k(k+\alpha+\beta+1))^rt^{2r}} \quad \mbox{ and } \quad b_k= \eta(atk), \quad \quad k=1,2,\ldots.
$$
And we note that $| \triangle^i \eta (atk) | \lesssim (at)^i$, for any $i$ a positive integer.\ By Lemma \ref{lem7} part a) and Lemma \ref{lem9}, respectively, we reach to
\begin{eqnarray*}
|\triangle^s \mu_{k,2} |  \lesssim \sum_{i=0}^s {s \choose i} \left(\frac{1}{k^{s+1}}+\frac{1}{k^{s}}\right) (kt)^s \lesssim \left(\frac{1}{k^{s+1}}+\frac{1}{k^{s}}\right).
\end{eqnarray*}
We still need to verify that the following holds true
\begin{eqnarray*}
\sup_j \left\{2^{j(s-1)} \sum_{k=2^j}^{2^{j+1}} |\bigtriangleup^s \mu_{k,2}|\right\} \lesssim c,
\end{eqnarray*}
for some constant $c$.
In fact, from previous estimates, we have
\begin{eqnarray*}
\sup_j \left\{2^{j(s-1)} \sum_{k=2^j}^{2^{j+1}} |\bigtriangleup^s \mu_{k,2}|\right\} \lesssim \sup_j 2^{j(s-1)} \left\{\sum_{k=2^j}^{2^{j+1}} \left(\frac{1}{k^{s+1}}+\frac{1}{k^{s}}\right) \right\}
			\lesssim  \sup_j \left\{\frac{1}{2^j} +1\right\}.
\end{eqnarray*}
Thus, Theorem \ref{Marcinkiewicz} implies that $\{\mu_{k,2}\}_k $ is a multiplier sequence and inequality (\ref{eq3}) holds.

Finally, we show that inequality (\ref{eq2}) holds from we showing that for $0<kt<\tau$ and $\tau>0$,  
\begin{equation*}
\mu_{k,3}:= \frac{(k(k+\alpha+\beta+1))^r t^{2r}}{1-m_r(k,t)} \, \eta(atk),\quad k=0,1,\ldots,
\end{equation*}
is a multiplier sequence fitting in Theorem \ref{Marcinkiewicz}.\ We observe that in Lemmas \ref{lem9} and \ref{lem7}, part a), we can replace sequence $\{a_k\}_k$ by $\{a_k^{-1}\}_k$.\ Also, Lemma \ref{lem7}, part b) fits in our context for the sequence bellow and
$$
a_k:=  \frac{1-m_r(k,t)}{(k(k+\alpha+\beta+1))^rt^{2r}} \geq c>0, k=0,1,\ldots,
$$ 
for some constant $c$ by Lemma \ref{eqmultiplier}.\ More than that  $a_k^{-1}\geq b^{-1}>0$, $k=0,1,\ldots$, where $b$ is the constant in formula (\ref{1}).\  Taking in account remark right after the statement of Lemma \ref{lem7} it is not hard to see $\{\mu_{k,3}\}_k$ is a multiplier sequence.\ The theorem is proved. \eop


\section{Application: decay of eigenvalues sequences}
\setcounter{equation}{0}

Our goal in this section is to prove both Theorem \ref{thmdecay2} and Corollary \ref{thmdecay1}.\ To present them we will first derive some additional technical results as following.\ We remind readers that the kernels $K$ we are dealing with satisfy all assumptions made Section 1.

We start noting that positivity of the kernel $K$assures that the operator ${\cal{L}}_K$ is positive and has a uniquely defined square root operator ${\cal{L}}_K^{1/2}$ whose kernel $K_{1/2}$ has the following series expansion
\begin{equation}\label{meiokernel}
K_{1/2}(x, y)=\sum_{k=0}^\infty \sum^{d^m_k}_{j=1} a^{1/2}_{k,j}\,  Y_{k,j}(x)\, Y_{k,j}(y), \quad  x,y\in \mathbb{M}.
\end{equation}
Both ${\cal{L}}_K$ and ${\cal{L}}_K^{1/2}$ are self-joint positive operators.\ The definition of the integral operator generated by $K$ makes easy to see that the spherical harmonics $Y_{k,j}$, $j=1,2, \ldots, d_k^m$ and $k=0,1,\ldots$, are all eigenvectors of the operator ${\cal{L}}_K$ associated to the eigenvalues $a_{k,j}$, respectively.\ Since we have made a monotonicity assumption on coefficients of $K$ it gives us an eigenvalue sequence ordering that is suitable for our analysis.\ 

For each $y \in \mathbb{M}$, the Fourier coefficients of the function $K^y:= K(\cdot, y)$ are $(\widehat{K^y})_{k,j}=a_{k,j}\,\overline{Y_{k,j}(y)}$, $j=1,2, \ldots, d_k^m$ and $k=0,1,\ldots$.\ Considering the kernel $K_{1/2}$ (formula (\ref{meiokernel})) in a similar way we have its Fourier coefficients $(\widehat{K_{1/2}^y})_{k,j}= a_{k,j}^{1/2}\,\overline{Y_{k,j}(y)}$,  $j=1,2, \ldots, d_k^m$ and $k=0,1,\ldots$, which implies that
\begin{equation}\label{dsum}
\int_{\mathbb{M}}s_k(K_{1/2}^y)\,dy=\sum_{j=1}^{d_k^m}a_{k,j}, \quad k=0,1,\ldots.
\end{equation}

The action of the fractional derivative on $K_{1/2}^y$ is given by
$$
\mathcal{B}^r(K_{1/2}^y) \sim \sum_{k=0}^{\infty}\sum_{j=1}^{d_k^m}a_{k,j}^{1/2}\,(k(k+\alpha+\beta+1))^{r/2}\,\overline{Y_{k,j}(y)}\,Y_{k,j}, \quad y\in \mathbb{M},
$$
and it permits us, with simple calculation, to derive the following $\left \|\mathcal{B}^r(K_{1/2}^y)\right\|_2^2 =\mathcal{B}^{2r,0}K(y,y)$, for any $y \in \mathbb{M}$.
Also, if $K$ is $(B,\beta)$-H\"{o}lder (see \cite{jordaosun}), then 
\begin{equation}\label{lemma2}
\int_{\mathbb{M}}\|S_t(K_{1/2}^y)-K_{1/2}^y\|_2^2\,dy\lesssim t^{\beta}, \quad y\in \mathbb{M}.
\end{equation}

\vspace{0,25cm}

\noindent \textbf{Proof of Theorem \ref{thmdecay2}} This proof can be found on the spherical setting in \cite{jordaosun}.\ Since from this point it is exactly the same one presented in this reference we just draw some steps of it.

By Proposition \ref{ditinequality1} for $p=q=2$ and $r=2$ and Theorem \ref{desejada} we have
$$\sum_{k=1}^{\infty}(\mbox{min}\{1, tk\})^{4}\,s_k(K_{1/2}^z)\lesssim \|S_t(K_{1/2}^z)-K_{1/2}^z\|_2^2, \quad z \in \mathbb{M},\,\, t\in(0,\pi).$$
Integrating both sides and making use of (\ref{lemma2}) we have $\sum_{k=0}^{\infty}(\mbox{min}\{1, tk\})^{4}\sum_{j=1}^{d_k^m}a_{k,j}\lesssim t^{\beta}$.\ Handling this inequality (for $t=1/n$) we get
$$
n^{\beta+m}\,a_n = n^{\beta+m-1}\sum_{k=n}^{2n-1} a_n \leq n^{\beta+m-1}\sum_{k=n}^{\infty}a_k\leq C_3, \quad n=1,2,\ldots,
$$
or, equivalently, $a_n =O(n^{-\beta-m})$, as $n\to \infty$.\ Returning to our original notation for the eigenvalues of $\mathcal{L}_K$
and recalling that $\{\lambda_{n}(\mathcal{L}_K)\}_n$ decreases to 0, we have that
$a_n = \lambda_{d_n^{m+1}}(\mathcal{L}_K)$, $n=1,2,\ldots$.\ In particular,
$$\lambda_{d_n^{m+1}}(\mathcal{L}_K) = O(n^{-\beta-m}), \quad n \to \infty.$$  Therefore, the decay in the statement of the theorem follows.\eop

\vspace{0,25cm}

\noindent \textbf{Proof of Corollary \ref{thmdecay1}} Most of steps in this proof are essentially repetitions of previous theorem that is why we omitted it here.\ By Proposition \ref{ditinequality1} ($p=q=2$ and the function $K_{1/2}^z$) we have
$$
\sum_{k=0}^{\infty}(\mbox{min}\{1, tk\})^{2r}s_k(K_{1/2}^z)\lesssim\left[\omega_{r}(K_{1/2}^z,t)_2\right]^2, \quad z\in \mathbb{M}, \,\, t \in (0,\pi).
$$
Since $K_{1/2}^z\in W^{2r}_2$, Proposition 4.2 in \cite{platonov1} asserts that $\omega_{r}(K_{1/2}^z,t)_2\lesssim\, t^{2r}\,\|\mathcal{B}^{r}(K_{1/2}^z)\|_2$, $z \in \mathbb{M}$.\ Then, we have
$$
\sum_{k=0}^{\infty}(\mbox{min}\{1, tk\})^{2r}\left(\int_{\mathbb{M}}s_k(K_{1/2}^z)\,dz\right)\lesssim t^{2r}
\int_{\mathbb{M}}\|\mathcal{B}^{r}(K_{1/2}^z)\|_2^2\,dz, \quad t \in (0,\pi).
$$
Since $B^{2r,0}K$ is the kernel of a trace-class operator $\|\mathcal{B}^{r}(K_{1/2}^z)\|_2^2$ is a nonnegative constant.\ Calculations analogous to before finishes the proof.\eop

\section{Kernel-based spaces and example}
\setcounter{equation}{0}
 
The last application to be presented corroborates with Theorem 6 in \cite{santin}.\ There the authors work on a context including just the Euclidean space and a bounded with smooth enough boundary domains and they recover decay rates for sequences of eigenvalues of integral operators from the decay of $n$-widths.\ Technique applied here is completely different from that since we make opposite way.

The\textit{ Kolmogorov $n$-width} of a subset $A$ of a Hilbert space $\left(H, \langle \cdot, \cdot\rangle_H\right)$ is defined as follows
\begin{equation}\label{kolmogorovn}
d_n(A;H):= \inf_{V_n\subset H} \,\, \sup_{f\in A} \,\, \inf_{f_n\in V_n} \,\, \|f-f_n\|_H,
\end{equation}
where $\|\cdot\|_H$ in the induced norm by the inner product in $H$ and the first infimum above is taken over all subspaces $V_n$ having dimension $n$ in $H$.

Additionally, under assumptions made here and continuity of the kernel $K$ the classical Mercer's Theorem assures that integral operator has a sequence of positive eigenvalues $\{\lambda_i\}$ ordered in a decreasing way, related to a sequence of eigenfunctions $\{\varphi_i\}$.\ More than that the kernel $K$ can be written as
\begin{equation}\label{mercer}
K(x,y)=\sum_{i=1}^{\infty}\lambda_i\varphi_i(x)\varphi_i(y), \quad x,y\in \mathbb{M},
\end{equation}
where the sum is absolutely and uniformly convergent.\ This result permits us to characterize $\left(\mathcal{H}_K, \langle \cdot, \cdot\rangle_K\right)$, which is the unique reproducing kernel Hilbert space (RKHS) attached to $K$, since $\{\sqrt{\lambda_i} \, \varphi_i\}$ is an orthonormal basis of it, (a complete reference for this basic facts are papers authored by R. Schaback).

Let us denote by $d_n$ the $n$-width $d_n(S(\mathcal{H}_K);L^2(\mathbb{M}))$.\ Then we have (\cite[Corollary 2.6]{pinkus})
\begin{equation}
d_n = \inf_{V_n\subset L^2} \,\, \sup_{f\in S(\mathcal{H}_K)} \,\, \|f-\mathcal{P}_n(f)\|_2= \sqrt{\lambda_{n+1}},
\end{equation}
where $S(\mathcal{H}_K)$ is unit ball in $\mathcal{H}_K$,  and the projections $\mathcal{P}_n$ are defined  by $\mathcal{P}_n(f):= \sum_{i=1}^{n}\langle f, h_i\rangle_K h_i$, $n=1,2,\ldots$,
where $\{h_i : i=1,2,\ldots, n\}$ is an orthonormal basis of $V_n$.
Moreover,  $H_n=\, \mbox{span}\,\{\sqrt{\lambda_i} \, \varphi_i : i=1,\ldots, n\}$, is the unique optimal space.

\ According to Santin and Schaback (\cite[p. 979]{santin}) if we consider a positive definite and symmetric kernel $K: \mathbb{M}\times \mathbb{M}\longrightarrow \mathbb{R}$, the Kolmogorov $n$-width $d_n$, defined in formula (\ref{kolmogorovn}) is equivalent to 
$$
\kappa_n:= \inf_{V_n\subset \mathcal{H}_K} \,\, \sup_{f\in S(\mathcal{H}_K)} \|f-\mathcal{P}_n(f)\|_2,
$$
where the infimum above is taken over all subspaces $V_n$ having dimension $n$ in $\mathcal{H}_K$.\ We bring up the following characterization for the Kolmogorov $n$-width, on this context.\ 

\begin{prop}\label{kolmogorovcarc}(\cite[Theorem 2]{santin}) If $H_n=\, \mbox{span}\,\{\sqrt{\lambda_i} \, \varphi_i : i=1,\ldots, n\}$, the subspace of $\mathcal{H}_K$, then $\kappa_n=\sqrt{\lambda_{n+1}}$.\ Moreover, $H_n$ is the unique optimal space.
\end{prop}

\begin{thm}\label{Kolmogorovn}
Let $K$ be a continuous, positive definite and symmetric kernel defined on $\mathbb{M}$.\ If $K$ satisfies the $(B,\beta)$ - H\"{o}lder condition, then
$$
\kappa_n = O((n+1)^{-1/2-\beta/2m}), \quad n\to\infty.
$$
\end{thm}

The proof from a simple application of Proposition \ref{kolmogorovcarc} and Theorem \ref{thmdecay2}.

\subsection{A concrete case: example} 
\setcounter{equation}{0}

The example bellow is a constructive way to consider a kernel to show the decay rates for the integral operator generating for it fits into assumptions of our theorems.\ 

Let $\epsilon>0$ be fixed and suppose $m\epsilon>1$ and $K$ is a kernel having expansion in the form

\begin{equation}\label{exemplo}
K(x,y)\sim  1+\sum_{n=1}^{\infty} \frac{c_n}{n^{m(1+\epsilon)+2r-1}} P_n^{(\alpha,\beta)}(\cos t), \quad x,y\in\mathbb{M},
\end{equation}
where $\cos t= d(x,y)$, and 
$$
c_n= \frac{\Gamma(\beta+1)(2n+\alpha+\beta+1)\Gamma(n+\alpha+\beta+1)}{\Gamma(\alpha+\beta+2)\Gamma(n+\beta+1)}.
$$
The  harmonic expansion of $K$ is easily obtained with the help of the addition formula and it is not hard to see that
\begin{eqnarray*}
1+ \sum_{n=1}^{\infty} \frac{c_n}{n^{m(1+\epsilon)+2r-1}}  &\leq & 1+ C\,\sum_{n=1}^{\infty} \frac{1}{n^{m\epsilon +2r}} <\infty, 
\end{eqnarray*}
for some constant $C$ depending on $m$.\ It means that the series expansion of $K$ (\ref{exemplo}) converges uniformly to $K(x, y)$ and then $K$ is continuous.\ Also its integral operator is positive, since $K$ is positive definite.\ The integral operator generated by  $\mathcal{B}^{2r,0}K$ is trace-class and $K$ fits into Corollary \ref{thmdecay1}.\ Therefore $\lambda_n({\cal{L}}_K)=O( n^{-1-2r/m})$, as $n\to \infty$.\ Also, Theorem \ref{Kolmogorovn} is applicable and then the Kolmogorov $n$-width of $H_n$ for example above decay as
$$
\kappa_n = O((n+1)^{-1/2-\beta/2m}), \quad n\to\infty.
$$
%
%
\vspace{-0.25cm}

\bibliographystyle{amsplain}

\begin{thebibliography}{99}
%
\bibitem{askey} R. Askey,  \textit{Orthogonal polynomials and special functions}. Society for Industrial and Applied Mathematics, Philadelphia, Pa., 1975.
%
\bibitem{berg} C. Berg, A.P. Peron, E. Porcu, \textit{Orthogonal expansions related to compact Gelfand pairs}. Expositiones Mathematicae, 2017. 
%
\bibitem{Bonami} A. Bonami, J.L. Clerc, \textit{Sommes de C\`esaro et multiplicateurs des d\'{e}veloppements en harmoniques sph\'{e}riques}. Trans. Amer. Math. Soc. \textbf{183} (1973), 223--263.
%
\bibitem{dai} G. Brown, F. Dai, \textit{Approximation of smooth functions on compact two-point homogeneous spaces}. J. Funct. Anal. \textbf{220} (2005), no. 2, 401--423.
%
\bibitem{mathcomp} M.H. Castro, V.A. Menegatto, \textit{Eigenvalue decay of positive integral operators on the sphere}. Math. Comp. \textbf{81} (2012), no. 280, 2303--2317.
%
\bibitem{daidi} F. Dai, Z. Ditzian, \textit{ Combinations of multivariate averages.}  J. Approx. Theory \textbf{131} (2004), no. 2, 268--283.
%
\bibitem{dai-xu} F. Dai, Y. Xu, \textit{ Approximation theory and harmonic analysis on spheres and balls.} Springer Monographs in Mathematics. Springer, New York, 2013.
%
\bibitem{dai-wang} F. Dai, K. Wang, C. Yu, \textit{On a conjecture of Ditzian and Runovskii}, J. Approx. Theory \textbf{118} no. 2 (2002) 202--224.
%
\bibitem{ditzian1} Z. Ditzian, \textit{Fractional derivatives and best approximation}. Acta Math. Hungar. \textbf{81} (1998), no. 4, 323--348.
%
\bibitem{ditzian2}  Z. Ditzian, \textit{Relating smoothness to expressions involving Fourier coefficients or to a Fourier transform}. J. Approx. Theory \textbf{164} (2012), no. 10, 1369--1389. 
%
\bibitem{jordaomen} T. Jord\~{a}o, V. A. Menegatto, \textit{Estimates for Fourier sums and eigenvalues of integral operators via multipliers on the sphere}. Proc. Amer. Math. Soc. \textbf{144} (2016), no. 1, 269--283. 
%
\bibitem{jordaosun} T. Jord\~{a}o, V. A. Menegatto, X. Sun, \textit{Eigenvalue sequences of positive integral operators and moduli of smoothness}. Springer, Cham, \textbf{83} (2014), p. 239--254.
%
\bibitem{kuhn} T. K\"{u}hn, \textit{Eigenvalues of integral operators with smooth positive definite kernels}.\ Arch. Math. (Basel) \textbf{49} (1987), no. 6, 525--534.
%
\bibitem{kushpel2} A. Kushpel, S.A. Tozoni, \textit{Entropy and widths of multiplier operators on two-point homogeneous spaces}. Constr. Approx., \textbf{35} (2012), no.2, 137--180.
%
\bibitem{pinkus} A. Pinkus, \textit{$n$-Widths in Approximation Theory}. Springer-Verlag, Berlin, 1985. 
%
\bibitem{platonov2}  S.S. Platonov, \textit{Approximations on compact symmetric spaces of rank 1}. (Russian) Mat. Sb. \textbf{188} (1997), no. 5, 113--130; translation in Sb. Math. \textbf{188} (1997), no. 5, 753--769.
%
\bibitem{platonov1} S.S. Platonov, \textit{Some problems in the theory of the approximation of functions on compact homogeneous manifolds.} (Russian) Mat. Sb. \textbf{200} (2009), no. 6, 67--108; translation in Sb. Math. \textbf{200} (2009), no. 5-6, 845--885.
%
\bibitem{santin} G. Santin, R. Schaback, \textit{Approximation of eigenfunctions in kernel-based spaces.}\ { Adv. Comput. Math.} \textbf{42} (2016), no. 4, 973--993.
%
\bibitem{schoenberg}I. J. Schoenberg, \textit{Positive definite functions on spheres}.\ {Duke Math. J.}  \textbf{9} (1942), 96--108.
%
\bibitem{szego} G. Szeg\"o,  \textit{Orthogonal Polynomials}. Amer. Math. Soc. New York, 1939.
%
\bibitem{wang} H-C. Wang, \textit{Two point homogeneous spaces}. Ann. Math. \textbf{ (2) 55}, (1952), 177--191.
\end{thebibliography}


\vspace*{0.25cm}

\noindent A. O. Carrijo\,\,\, \& \,\,\,T. Jord\~{a}o\\
Departamento de
Matem\'atica,\\
ICMC - Universidade de S\~{a}o Paulo, \\
13566-590 - S\~ao Carlos - SP, Brazil.

\end{document}